# Octahedral Galois representations arising from **Q**-curves of degree 2
J. Fernández, J-C. Lario and A. Rio

June 16, 2001

**Abstract**

Generically, one can attach to a **Q**-curve $C$ octahedral representations $\rho\colon \mathrm{Gal}(\overline{\mathbf{Q}}/\mathbf{Q}) \longrightarrow \mathrm{GL}_2(\overline{\mathbf{F}}_3)$ coming from the Galois action on the 3-torsion of those abelian varieties of $\mathrm{GL}_2$-type whose building block is $C$. When $C$ is defined over a quadratic field and has an isogeny of degree 2 to its Galois conjugate, there exist such representations $\rho$ having image into $\mathrm{GL}_2(\mathbf{F}_9)$. Going the other way, we can ask which mod 3 octahedral representations $\rho$ of $\mathrm{Gal}(\overline{\mathbf{Q}}/\mathbf{Q})$ arise from **Q**-curves in the above sense. We characterize those arising from quadratic **Q**-curves of degree 2. The approach makes use of Galois embedding techniques in $\mathrm{GL}_2(\mathbf{F}_9)$, and the characterization can be given in terms of a quartic polynomial defining the $\mathcal{S}_4$-extension of **Q** corresponding to the projective representation $\overline{\rho}$.

## 1 Introduction

The motivation of the present paper is nailed in the context of a conjecture of Serre [Ser87]: odd irreducible residual 2-dimensional linear representations of $\mathrm{Gal}(\overline{\mathbf{Q}}/\mathbf{Q})$ should arise from modular abelian varieties. Rephrasing it, in a weaker form, we can say that they should arise from abelian varieties of $\mathrm{GL}_2$-type. Our aim is to explore the minimal dimensions for these abelian varieties, as well as for their building blocks. Here, we will restrict ourselves to certain mod 3 octahedral Galois representations. For previous work in this direction, we refer to [SBT97] and [LR95]. Before explaining the plan of the paper, we begin by reviewing the basic ingredients that will be used.

• **Q-curves and abelian varieties of** $\mathrm{GL}_2$**-type.** Throughout, the term **Q**-*curve* will stand for a non-CM elliptic curve defined over $\overline{\mathbf{Q}}$ which is isogenous to all its Galois conjugates. An isogeny between a **Q**-curve and one of its conjugates is then determined up to multiplication by rational numbers; we say that the **Q**-curve is *completely defined* over a Galois number field if all its conjugates and the isogenies between them are defined over that field. By a **Q**-*curve of degree* $N$ we will mean a **Q**-curve defined over a quadratic field, with non-rational $j$-invariant, and having an isogeny of degree $N$ to its conjugate curve.

An abelian variety $A$ defined over $\mathbf{Q}$ is said to be *of* $\mathrm{GL}_2$-*type* if its algebra of $\mathbf{Q}$-endomorphisms $\mathbf{Q} \otimes \mathrm{End}_{\mathbf{Q}}(A)$ is a number field of degree $\dim(A)$. Such an abelian variety is $\mathbf{Q}$-simple and it is isogenous to the power of some absolutely simple abelian variety, which is called its *building block* (see [Rib94] and [Pyl95]). In [Rib92], Ribet characterizes $\mathbf{Q}$-curves as the elliptic curves $C$ over $\overline{\mathbf{Q}}$ which are the building block of some abelian variety $A$ of $\mathrm{GL}_2$-type; we refer to this situation by saying that $A$ is an abelian variety *attached to* the $\mathbf{Q}$-curve $C$. Conjecturally, $\mathbf{Q}$-curves are the (non-CM) elliptic curves over $\overline{\mathbf{Q}}$ which are quotient of some modular curve $X_1(N)$.

• **Octahedral representations.** Let $k$ be a perfect field. Given a continuous representation $\rho\colon \mathrm{Gal}(\overline{k}/k) \to \mathrm{GL}(V)$, let $\overline{\rho}$ denote its associated projective representation; i.e., the composition of $\rho$ with the canonical map $\pi\colon \mathrm{GL}(V) \to \mathrm{PGL}(V)$. Usually, we will denote by $K_\rho$ and $K_{\overline{\rho}}$ the fixed fields of $\rho$ and $\overline{\rho}$, respectively. Such a representation $\rho$ is called *octahedral* if the image of $\overline{\rho}$ is isomorphic to the symmetric group $\mathcal{S}_4$. Mostly, we will be dealing with odd representations $\rho\colon \mathrm{Gal}(\overline{\mathbf{Q}}/\mathbf{Q}) \to \mathrm{GL}_2(\mathbf{F}_9)$ such that $\overline{\rho}\colon \mathrm{Gal}(\overline{\mathbf{Q}}/\mathbf{Q}) \to \mathrm{PGL}_2(\mathbf{F}_9)$ has image isomorphic to $\mathrm{PGL}_2(\mathbf{F}_3) \simeq \mathcal{S}_4$.

• **Principal quartics.** Let $K_1/k$ be an extension of number fields of degree 4. We will say that $K_1/k$ is a *principal quartic* if there is a polynomial $f(X) = X^4 + bX + c \in k[X]$ with splitting field the Galois closure of $K_1$ in $\overline{k}$. We will also say that $f$ is a *principal* quartic polynomial. Whenever $\mathrm{Gal}(f)$ is isomorphic to $\mathcal{S}_4$, the mod 3 octahedral representations attached to $f$ will also be called *principal*.

In section 2 we explain how to attach octahedral Galois representations to $\mathbf{Q}$-curves. Section 3 deals with octahedral embedding problems in $\mathrm{GL}_2(\mathbf{F}_9)$, and section 4 is devoted to linking them to principal quartics. The characterization of octahedral extensions of $\mathbf{Q}$ arising from $\mathbf{Q}$-curves of degree 2 and the $\mathbf{Q}$-endomorphism algebras of the abelian varieties of $\mathrm{GL}_2$-type attached to them are explained in sections 5 and 6. Finally, section 7 contains several numerical examples to illustrate the different embedding problems considered in the previous sections.

## 2 Octahedral representations attached to Q-curves

Let $C$ be a $\mathbf{Q}$-curve defined over a number field $k$ and let $p$ be a prime number. Set $\mathrm{G}_k = \mathrm{Gal}(\overline{k}/k)$ and consider the representation

$$\rho_{C,p}\colon \mathrm{G}_k \to \mathrm{Aut}(C[p]) \simeq \mathrm{GL}_2(\mathbf{F}_p)$$

given by the action on the $p$-torsion module of $C$. In this section we briefly explain, following [ES00] and [Rib92], a procedure to obtain continuous representations $G_{\mathbf{Q}} \longrightarrow \mathrm{GL}_2(\overline{\mathbf{F}}_p)$ whose restriction to $G_k$ is isomorphic to a twist of $\rho_{C,p}$. For $p = 3$, and assuming $\rho_{C,3}$ to be surjective, this procedure gives rise to octahedral representations.

Fix $\{\mu_\sigma\}_{\sigma \in G_{\mathbf{Q}}}$ a locally constant system of isogenies. That is, for every $\sigma \in G_{\mathbf{Q}}$ we choose an isogeny $\mu_\sigma \colon {}^\sigma C \to C$ in such a way that $\mu_\sigma = \mu_\tau$ whenever $\sigma$ and $\tau$ restrict to the same embedding of $k$ into $\overline{\mathbf{Q}}$. Also, we take $\mu_\sigma$ to be the identity for all $\sigma$ in $G_k$. The map

$$c \colon (\sigma, \tau) \longmapsto \mu_\sigma \, {}^\sigma\mu_\tau \, \mu_{\sigma\tau}^{-1}$$

is a 2-cocycle of the group $G_{\mathbf{Q}}$ with values in $\overline{\mathbf{Q}}^*$. If we let $G_{\mathbf{Q}}$ act trivially in $\overline{\mathbf{Q}}^*$, the corresponding cohomology group $H^2(G_{\mathbf{Q}}, \overline{\mathbf{Q}}^*)$ is trivial (see [Ser77] and [Que95]), so that there exists a continuous map $\alpha \colon G_{\mathbf{Q}} \to \overline{\mathbf{Q}}^*$ such that

$$c(\sigma, \tau) = \alpha(\sigma)\,\alpha(\tau)\,\alpha(\sigma\tau)^{-1}$$

for all $\sigma, \tau \in G_{\mathbf{Q}}$. In particular, the restriction of $\alpha$ to $G_k$ is a Galois character $\psi$. Let $E$ be the number field generated by the values of $\alpha$. For each prime $\mathfrak{p}$ of $E$ over $p$, we let $E_\mathfrak{p}$ be the completion of $E$ at $\mathfrak{p}$, viewed as a finite extension of $\mathbf{Q}_p$.

Consider now the $p$-adic Tate module $T_p(C)$ and the usual representation $\phi_{C,p}$ of $G_k$ attached to it. For each $\mathfrak{p}$ as above, $T_\mathfrak{p}(C) := E_\mathfrak{p} \otimes_{\mathbf{Z}_p} T_p(C)$ is then a 2-dimensional $E_\mathfrak{p}$-vector space with an action of $G_{\mathbf{Q}}$ given by

$$^\sigma(1 \otimes P) := \alpha(\sigma)^{-1} \otimes \mu_\sigma({}^\sigma P),$$

where $P$ denotes an element of $T_p(C)$. Thus, we have a representation

$$\phi_{C,\mathfrak{p}} \colon G_{\mathbf{Q}} \longrightarrow \mathrm{GL}_2(E_\mathfrak{p}).$$

By choosing appropriate bases for $T_p(C)$ and $T_\mathfrak{p}(C)$, we can assume that $\phi_{C,\mathfrak{p}}$ has image inside the subgroup $E^* \mathrm{GL}_2(\mathbf{Z}_p)$ and also that its restriction to $G_k$ equals $\psi^{-1} \otimes \phi_{C,p}$. Notice that $\phi_{C,\mathfrak{p}}$ depends not only on the prime $\mathfrak{p}$, but also on the splitting map $\alpha$. Observe also that any other splitting map differs from $\alpha$ by a character of $G_{\mathbf{Q}}$.

Associated with this splitting map, in [Rib92] Ribet attaches to $C$ an abelian variety $A$ of $\mathrm{GL}_2$-type having the number field $E$ as algebra of $\mathbf{Q}$-endomorphisms. It is constructed as a $\mathbf{Q}$-simple factor of the abelian variety $\mathrm{Res}_{L/\mathbf{Q}}(C)$ obtained by restriction of scalars from a certain Galois extension $L$ of $\mathbf{Q}$, depending on $\alpha$, over which the $\mathbf{Q}$-curve is completely defined. Moreover, in [Pyl95] it is shown that any abelian variety attached to $C$ is obtained, up to isogeny, from this construction applied to some splitting map. The Galois action on the $p$-adic Tate module $V_p(A)$ defines, for each prime $\mathfrak{p}$ of $E$ over $p$, a representation

$$\phi_{A,\mathfrak{p}} \colon G_{\mathbf{Q}} \longrightarrow \mathrm{GL}_2(E_\mathfrak{p})$$

whose isomorphism class has an attached reduction $\rho_{A,\mathfrak{p}}$ over the residue field $\mathbf{F}_\mathfrak{p}$ of $\mathfrak{p}$. For an account on $\mathfrak{p}$-adic representations attached to abelian varieties of $\mathrm{GL}_2$-type and their reductions mod $\mathfrak{p}$ we refer to [Rib92] and, in a more general setting, to [Rib76].

According to [ES00], $\phi_{C,\mathfrak{p}}$ is an odd continuous representation isomorphic to Ribet's $\phi_{A,\mathfrak{p}}$. Also, it has some conjugate over $\overline{\mathbf{Q}}_p$ with a reduction taking values into $\overline{\mathbf{F}}_p^* \mathrm{GL}_2(\mathbf{F}_p)$. From all the above, we have the following result:

**Proposition 2.1** *Let $C$ be a $\mathbf{Q}$-curve over a number field $k$, and $A$ over $\mathbf{Q}$ be an abelian variety of $\mathrm{GL}_2$-type associated with a splitting map $\alpha$ as above. For every prime ideal $\mathfrak{p}$ of $\mathbf{Q} \otimes \mathrm{End}_{\mathbf{Q}}(A)$ over a prime $p$, we have, attached to $C$ and $p$, an odd representation*

$$\rho_{A,\mathfrak{p}}\colon \mathrm{G}_{\mathbf{Q}} \longrightarrow \mathrm{GL}_2(\mathbf{F}_{\mathfrak{p}})$$

*such that:*

- *the image of $\overline{\rho}_{A,\mathfrak{p}}$ is conjugate inside $\mathrm{PGL}_2(\overline{\mathbf{F}}_p)$ to a subgroup of $\mathrm{PGL}_2(\mathbf{F}_p)$;*
- *its restriction $\rho_{A,\mathfrak{p}|k}$ to $\mathrm{G}_k$ is isomorphic to $\widetilde{\psi}^{-1} \otimes \rho_{C,p}$, where $\widetilde{\psi}$ is the reduction mod $\mathfrak{p}$ of $\psi = \alpha_{|k}$. In particular, $\overline{\rho}_{A,\mathfrak{p}|k}$ and $\overline{\rho}_{C,p}$ have the same fixed field.*

We will say that $\rho\colon \mathrm{G}_{\mathbf{Q}} \to \mathrm{GL}_2(\overline{\mathbf{F}}_p)$ *arises from a $\mathbf{Q}$-curve* if $\rho$ is isomorphic to a representation $\rho_{A,\mathfrak{p}}$ attached to a $\mathbf{Q}$-curve as described in the proposition. Notice that this property is invariant by twists; so we will also say that $\overline{\rho}$ arises from a $\mathbf{Q}$-curve and that $A$ is an abelian variety attached to $\overline{\rho}$.

The rest of the paper deals with the particular case in which $p$ is 3. Assume that $\rho_{C,3}$ is surjective or, equivalently, that the fixed field of $\overline{\rho}_{C,3}$ is an $\mathcal{S}_4$-extension of $k$. We recall that this field is obtained by adjoining to $k$ the $x$-coordinates of the 3-torsion points of $C$ using any Weierstrass equation for $C$ over $k$, and it contains the cubic roots of unity. With that starting point, any representation

$$\rho\colon \mathrm{G}_{\mathbf{Q}} \longrightarrow \mathrm{GL}_2(\overline{\mathbf{F}}_3)$$

attached to $C$ and 3 like in proposition 2.1 is octahedral. If $k$ is a quadratic field, then for any such $\rho$ the fixed field of $\overline{\rho}$ is the unique $\mathcal{S}_4$-extension of $\mathbf{Q}$ whose compositum with the quadratic $k$ equals the fixed field of $\overline{\rho}_{C,3}$.

Here we will only be concerned with $\mathbf{Q}$-curves over quadratic fields. It is a fact that any $\mathbf{Q}$-curve $C$ of degree 2 has always an attached abelian variety with $\mathbf{Q}$-endomorphism algebra $\mathbf{Q}\left(\sqrt{2}, \sqrt{-2}\right)$ (see [Que00] and section 6). In this number field there are two primes over 3, each one with residual degree 2. Hence, in the surjective case, there exist odd octahedral representations

$$\rho\colon \mathrm{G}_{\mathbf{Q}} \longrightarrow \mathrm{GL}_2(\mathbf{F}_9)$$

attached to $\mathbf{Q}$-curves of degree 2. Our main goal is to characterize octahedral Galois representations into $\mathrm{GL}_2(\mathbf{F}_9)$ arising from $\mathbf{Q}$-curves of degree 2.

# 3 Octahedral embedding problems in $\mathrm{GL}_2(\mathbf{F}_9)$

The aim of this section is to determine the possible images of odd octahedral representations $\mathrm{G}_\mathbf{Q} \to \mathrm{GL}_2(\mathbf{F}_9)$ and to relate them to certain Galois embedding problems given by central extensions of $\mathcal{S}_4$ with kernel a cyclic group of order a power of 2.

Let $K/\mathbf{Q}$ be an $\mathcal{S}_4$-extension, and consider a projective representation

$$\varrho\colon \mathrm{G}_\mathbf{Q} \twoheadrightarrow \mathrm{Gal}(K/\mathbf{Q}) \simeq \mathrm{PGL}_2(\mathbf{F}_3) \hookrightarrow \mathrm{PGL}_2(\overline{\mathbf{F}}_3).$$

A *lifting* of $\varrho$ is a representation $\rho\colon \mathrm{G}_\mathbf{Q} \to \mathrm{GL}_2(\overline{\mathbf{F}}_3)$ such that $\overline{\rho} = \varrho$, where $\overline{\rho}$ denotes as above the composition of $\rho$ with the canonical map $\pi\colon \mathrm{GL}_2(\overline{\mathbf{F}}_3) \to \mathrm{PGL}_2(\overline{\mathbf{F}}_3)$. If we denote by $K_\rho$ the fixed field of $\rho$, we have a commutative diagram

$$\begin{array}{ccccccccc}
1 & \longrightarrow & \mathrm{G}_K & \longrightarrow & \mathrm{G}_\mathbf{Q} & & & & \\
& & \downarrow & & \downarrow & \searrow & & & \\
1 & \longrightarrow & \mathrm{Gal}(K_\rho/K) & \longrightarrow & \mathrm{Gal}(K_\rho/\mathbf{Q}) & \longrightarrow & \mathrm{Gal}(K/\mathbf{Q}) & \longrightarrow & 1 \\
& & \downarrow & & \downarrow & & \downarrow & & \\
1 & \longrightarrow & \overline{\mathbf{F}}_3^* & \longrightarrow & \mathrm{GL}_2(\overline{\mathbf{F}}_3) & \xrightarrow{\pi} & \mathrm{PGL}_2(\overline{\mathbf{F}}_3) & \longrightarrow & 1.
\end{array}$$

The degree of the cyclic extension $K_\rho/K$ is called the *index* of the lifting $\rho$. By a result of Tate (see [Ser77] and [Que95]), such a lifting $\rho$ always exists. All liftings of $\varrho$ are then $\rho \otimes \varepsilon$, for all characters $\varepsilon\colon \mathrm{G}_\mathbf{Q} \to \overline{\mathbf{F}}_3^*$. In particular, we can define $\varrho$ to be *odd* if $\det \rho$ is an odd character. This is equivalent to saying that the fixed field of $\det \rho$ is a CM field, and also to saying that the fixed field $K$ of $\varrho$ is not real.

The results of [Que95] show that $\varrho$ has always a lifting $\rho$ of index $2^r$ for some positive integer $r$. Then, $\mathrm{Im}\rho$ is a central extension of $\mathrm{Im}\overline{\rho} \simeq \mathcal{S}_4$ by a cyclic group $C_{2^r}$ of order $2^r$. As an element of the cohomology group $H^2(\mathcal{S}_4, C_{2^r})$ this image has only two possibilities, denoted by $2^r\mathcal{S}_4^+$ and $2^r\mathcal{S}_4^-$. Besides, up to twists of $\rho$ preserving the fixed field $K_\rho$, we can assume $\det \rho$ to have order either $2^r$ or $2^{r-1}$, depending on the quadratic subextension of the cyclic extension $K_{\det \rho}/\mathbf{Q}$ (see remark 3.1 below).

Before reproducing more precisely (in propositions 3.1 and 3.2 below) the relationship between liftings of $\varrho$ and solutions to the embedding problems associated with the above extensions of $\mathcal{S}_4$, we need some more cohomology.

For any Galois character $\varepsilon\colon \mathrm{G}_\mathbf{Q} \to F^*$, with $F$ an algebraically closed field of characteristic different from 2, let $[\varepsilon]$ be the pull-back by $\varepsilon$ of the exact sequence

$$\begin{array}{ccccccccc}
1 & \longrightarrow & \{\pm 1\} & \longrightarrow & F^* & \longrightarrow & F^* & \longrightarrow & 1. \\
& & & & x & \longmapsto & x^2 & &
\end{array}$$

It is the element of $H^2(G_{\mathbf{Q}}, \{\pm 1\}) \simeq \mathrm{Br}_2(\mathbf{Q})$ giving the obstruction to the existence of a character $\psi \colon G_{\mathbf{Q}} \to F^*$ such that $\psi^2 = \varepsilon$. It can also be interpreted as the obstruction to embedding the fixed field $K_\varepsilon$ of $\varepsilon$ into a cyclic extension of $\mathbf{Q}$ having twice its degree.

Let now $[s_4^-]$ be the pull-back by $\varrho$ of the exact sequence

$$1 \longrightarrow \{\pm 1\} \longrightarrow \mathrm{SL}_2(\overline{\mathbf{F}}_3) \longrightarrow \mathrm{PSL}_2(\overline{\mathbf{F}}_3) = \mathrm{PGL}_2(\overline{\mathbf{F}}_3) \longrightarrow 1.$$

It is the element of $H^2(G_{\mathbf{Q}}, \{\pm 1\}) \simeq \mathrm{Br}_2(\mathbf{Q})$ giving the obstruction to the existence of a lifting of $\varrho$ with trivial determinant. The component at $\infty$ of $[s_4^-]$ is $-1$ if and only if $\varrho$ is odd. This can be deduced, for instance, from the following result of Tate ([Ser77], Theorem 6).

**Theorem 3.1** *Given a character $\varepsilon \colon G_{\mathbf{Q}} \to \overline{\mathbf{F}}_3^*$, the projective representation $\varrho$ has a lifting with determinant $\varepsilon$ if and only if $[\varepsilon] = [s_4^-]$.*

Let $K_1/\mathbf{Q}$ be a quartic extension with normal closure $K/\mathbf{Q}$ and discriminant $d_{K_1}$, and let $w \in \mathrm{Br}_2(\mathbf{Q})$ be the Witt invariant of the quadratic form $\mathrm{Tr}_{K_1/\mathbf{Q}}(x^2)$. We have then the formula $[s_4^-] = w \otimes (-2, d_{K_1})$, due to Serre and Vila (cf. [Ser84], [Vil88]). Let also $\varepsilon_M \colon G_{\mathbf{Q}} \to \overline{\mathbf{F}}_3^*$ be the Galois character attached to $M = \mathbf{Q}(\sqrt{d_{K_1}})$, the only quadratic subfield of $K$. The next propositions come straightforward from the above theorem, Serre-Vila's formula and the results of [Que95].

**Proposition 3.1** *The embedding problem*

$$2^r \mathcal{S}_4^+ \longrightarrow \mathcal{S}_4 \simeq \mathrm{Gal}(K/\mathbf{Q})$$

*is solvable if and only if there exists a character $\varepsilon \colon G_{\mathbf{Q}} \to \overline{\mathbf{F}}_3^*$ such that*

$(a_+)$ $\varepsilon^{2^{r-1}} = \varepsilon_M$.

$(b_+)$ $[\varepsilon] = w \otimes (-2, d_{K_1})$.

*Moreover, $\rho \mapsto K_\rho$ (resp. $\rho \mapsto \det \rho$) is a surjective map from the set of liftings of $\varrho$ of index $2^r$ and determinant of order $2^r$ to the set of solutions to the embedding problem (resp. to the set of Galois characters $\varepsilon$ satisfying conditions $(a_+)$ and $(b_+)$).*

**Proposition 3.2** *The embedding problem*

$$2^r \mathcal{S}_4^- \longrightarrow \mathcal{S}_4 \simeq \mathrm{Gal}(K/\mathbf{Q})$$

*is solvable if and only if there exists a character $\varepsilon \colon G_{\mathbf{Q}} \to \overline{\mathbf{F}}_3^*$ such that*

$(a_-)$ $\varepsilon = 1$ for $r = 1$, $\varepsilon^{2^{r-2}}$ *is a quadratic character different from $\varepsilon_M$ for $r \geq 2$.*

$(b_-)$ $[\varepsilon] = w \otimes (-2, d_{K_1})$.

*Moreover, $\rho \mapsto K_\rho$ (resp. $\rho \mapsto \det \rho$) is a surjective map from the set of liftings of $\varrho$ of index $2^r$ and determinant of order $2^{r-1}$ to the set of solutions to the embedding problem (resp. to the set of Galois characters $\varepsilon$ satisfying conditions $(a_-)$ and $(b_-)$).*

Observe that $[s_4^-]$ is also the obstruction to the solvability of the embedding problem $2\mathcal{S}_4^- \to \mathcal{S}_4 \simeq \mathrm{Gal}(K/\mathbf{Q})$. In particular, this problem is never solvable when $K$ is not real.

**Remark 3.1** For a lifting $\rho$ as in the propositions, the following diagrams show where the fixed field of $\varepsilon = \det \rho$ lies in each case. Note that the extension $K_\rho/K_{\det \rho}$ has Galois group $\mathrm{SL}_2(\overline{\mathbf{F}}_3) \cap \mathrm{Im}\rho$, with order $|\mathcal{S}_4|$ in the first case and $2|\mathcal{S}_4|$ in the second case.

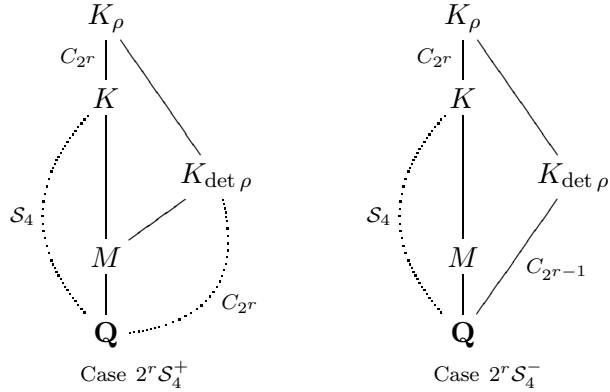

Case $2^r\mathcal{S}_4^+$     Case $2^r\mathcal{S}_4^-$

All liftings of $\varrho$ of index $2^r$ and determinant of order $2^r$ (resp. order $2^{r-1}$) are then $\rho \otimes \psi$, for all characters $\psi$ of order dividing $2^r$ (resp. of order dividing $2^r$ and, for $r > 1$, such that $\psi^{2^{r-1}}$ is different from the quadratic characters $(\det \rho)^{2^{r-2}}$ and $\varepsilon_M(\det \rho)^{2^{r-2}}$).

We are interested in liftings of $\varrho$ into $\mathrm{GL}_2(\mathbf{F}_9)$. The restriction of such a lifting $\rho$ to $G_K$ takes values in $\mathbf{F}_9^*$, so that $\mathrm{Im}\rho \simeq 2^r\mathcal{S}_4^+$ or $\mathrm{Im}\rho \simeq 2^r\mathcal{S}_4^-$ for $r$ either 1, 2 or 3; we exclude the case $r = 0$ because $\varrho$ cannot have trivial liftings, as one can easily see from the character table of $\mathcal{S}_4$. The actual possibilities for $\mathrm{Im}\rho$ are computed in the next proposition. We must first introduce some notation relative to the groups $\mathrm{GL}_2(\mathbf{F}_3)$ and $\mathbf{F}_9^*$. The first one is the subgroup $\langle S, T \rangle$ of $\mathrm{GL}_2(\mathbf{F}_9)$ generated by the matrices

$$S := \begin{pmatrix} 1 & 0 \\ 0 & -1 \end{pmatrix} \quad \text{and} \quad T := \begin{pmatrix} 1 & -1 \\ 1 & 0 \end{pmatrix}.$$

As for the multiplicative group of $\mathbf{F}_9$, if we let $\zeta$ and $i$ be respectively an eighth and a fourth primitive root of unity in $\overline{\mathbf{F}}_3$, then its non-trivial subgroups are $\langle \zeta \rangle = \mathbf{F}_9^*$, $\langle \zeta^2 \rangle = \langle i \rangle$ and $\langle \zeta^4 \rangle = \langle -1 \rangle = \mathbf{F}_3^*$. Notice also that it makes sense to define $\det \colon \mathrm{PGL}_2(\mathbf{F}_3) \to \{\pm 1\}$ as the natural morphism induced by the determinant $\det \colon \mathrm{GL}_2(\mathbf{F}_3) \to \{\pm 1\}$.

**Proposition 3.3** *The image of an octahedral representation*

$$\rho\colon G_{\mathbf{Q}} \longrightarrow \mathrm{GL}_2(\mathbf{F}_9)$$

*is conjugate to one of the following subgroups of* $\mathbf{F}_9^* \,\mathrm{GL}_2(\mathbf{F}_3)$:

$$\begin{aligned}
\mathrm{GL}_2(\mathbf{F}_3) = & \;\langle S, T\rangle \simeq 2\,\mathcal{S}_4^+,\\
& \;\langle \zeta\,S, T\rangle \simeq 4\,\mathcal{S}_4^+,\\
& \;\langle i\,S, T\rangle \simeq 2\,\mathcal{S}_4^-,\\
\langle i\rangle\,\mathrm{GL}_2(\mathbf{F}_3) = & \;\langle S, i\,T\rangle \simeq 4\,\mathcal{S}_4^-,\\
\mathbf{F}_9^*\,\mathrm{GL}_2(\mathbf{F}_3) = & \;\langle S, \zeta\,T\rangle \simeq 8\,\mathcal{S}_4^-.
\end{aligned}$$

*The third case does never occur if $\rho$ is odd. The outer automorphism groups of the five groups above are:*

$$\begin{aligned}
\mathrm{Out}(2\,\mathcal{S}_4^+) &\simeq \langle f_1\rangle,\\
\mathrm{Out}(4\,\mathcal{S}_4^+) &\simeq \langle f_1\rangle \times \langle f_2\rangle,\\
\mathrm{Out}(2\,\mathcal{S}_4^-) &\simeq \langle \varphi\rangle,\\
\mathrm{Out}(4\,\mathcal{S}_4^-) &\simeq \langle \varphi\rangle \times \langle f_1\rangle,\\
\mathrm{Out}(8\,\mathcal{S}_4^-) &\simeq \langle \varphi\rangle \times \langle f_1\rangle \times \langle f_2\rangle,
\end{aligned}$$

*where $\varphi$, $f_1$ and $f_2$ are involutions respectively defined by*

$$\begin{aligned}
M &\longmapsto \det(\pi(M))\,M\,,\\
M &\longmapsto \det(M)\,M\,,\\
M &\longmapsto \det(M)^2 M\,.
\end{aligned}$$

*Proof.* There is only one conjugacy class of subgroups isomorphic to $\mathcal{S}_4$ in $\mathrm{PGL}_2(\mathbf{F}_9)$. From this fact, together with the observation that two subgroups $H_1$ and $H_2$ of $\mathrm{GL}_2(\mathbf{F}_9)$ have the same image by $\pi$ if and only if $\mathbf{F}_9^*\,H_1$ equals $\mathbf{F}_9^*\,H_2$, one deduces that $\mathrm{Im}\rho$ is conjugate in $\mathrm{GL}_2(\mathbf{F}_9)$ to some subgroup of $\mathbf{F}_9^*\,\mathrm{GL}_2(\mathbf{F}_3)$. The projection of this last group to $\mathrm{PGL}_2(\overline{\mathbf{F}}_3)$ is the canonical embedding of $\mathrm{PGL}_2(\mathbf{F}_3)$, so we can assume that $\mathrm{Im}\rho$ has the same image by $\pi$ as $\mathrm{GL}_2(\mathbf{F}_3)$, and we must see then that it equals one of the groups in the statement.

Under this assumption, $\mathrm{Im}\rho$ contains the subgroup $H_{j,k} := \langle \zeta^j S, \zeta^k T\rangle$ for some $j, k \in \{0, 1, \ldots, 7\}$. Since the orders of $S$ and $T$ are 2 and 6 respectively, we have that some (even) power of the matrix $\zeta^j S$ becomes the homothety $-1$, and the same happens to the matrix $\zeta^k T$, except for the cases $j, k \in \{0, 4\}$, which correspond to the group $\mathrm{GL}_2(\mathbf{F}_3)$; hence, we always have $\mathbf{F}_3^* \subset H_{j,k}$. Using also the relations $(\zeta\,S)^3 = \zeta^3 S$, $(\zeta^3 S)^3 = \zeta\,S$, $(\zeta\,T)^7 = -\zeta^3 T$ and $(\zeta^3 T)^7 = -\zeta\,T$, we can restrict ourselves to the nine cases $j, k \in \{0, 1, 2\}$, which in fact reduce to five:

$$\begin{aligned}
G_1 &:= H_{0,0} = \mathrm{GL}_2(\mathbf{F}_3),\\
G_2 &:= H_{1,0} = H_{1,2},\\
G_3 &:= H_{2,0},\\
G_4 &:= H_{0,2} = H_{2,2} = \langle i\rangle\,\mathrm{GL}_2(\mathbf{F}_3),\\
G_5 &:= H_{0,1} = H_{1,1} = H_{2,1} = \mathbf{F}_9^*\,\mathrm{GL}_2(\mathbf{F}_3).
\end{aligned}$$

Since $G_j \subset \mathrm{Im}\rho \subset G_5$ for some $j$ and $\pi(G_j) = \pi\left(\mathrm{GL}_2(\mathbf{F}_3)\right) = \pi(\mathrm{Im}\rho)$, the only possibilities for $\mathrm{Im}\rho$ are $G_j$, $\langle i \rangle G_j$ or $\mathbf{F}_9^* G_j$. The subgroups $G_2$ and $G_4$ are maximal in $G_5$. Moreover, one can easily see that $\langle i \rangle G_1 = \langle i \rangle G_3 = G_4$ and $\mathbf{F}_9^* G_1 = \mathbf{F}_9^* G_3 = G_5$. Hence, $\mathrm{Im}\rho$ equals one of the five groups above.

Each $G_j$ corresponds to an element of $H^2(\mathcal{S}_4, C_{2^r})$ for some $r = 1, 2, 3$ and can be realized as a Galois group over $\mathbf{Q}$ (see [Son91]), so that $G_j \simeq 2^r \mathcal{S}_4^+$ or $G_j \simeq 2^r \mathcal{S}_4^-$. The isomorphism class of $G_j$ is determined by the orders of $G_j$ and $G_j \cap \mathrm{SL}_2(\mathbf{F}_9)$, which are displayed in the following table.

| $j$ | $|G_j|$ | $|G_j \cap \mathrm{SL}_2(\mathbf{F}_9)|$ |
|---|---|---|
| 1 | $2\,|\mathcal{S}_4|$ | $|\mathcal{S}_4|$ |
| 2 | $4\,|\mathcal{S}_4|$ | $|\mathcal{S}_4|$ |
| 3 | $2\,|\mathcal{S}_4|$ | $2\,|\mathcal{S}_4|$ |
| 4 | $4\,|\mathcal{S}_4|$ | $2\,|\mathcal{S}_4|$ |
| 5 | $8\,|\mathcal{S}_4|$ | $2\,|\mathcal{S}_4|$ |

Notice that $G_1 \cap \mathrm{SL}_2(\mathbf{F}_9) = \mathrm{SL}_2(\mathbf{F}_3)$ and that $G_3 \subset \mathrm{SL}_2(\mathbf{F}_9)$. The only non-obvious row is the second one, for which some computations using GAP have been made. The rest of the proposition is an easy calculation. $\square$

**Remark 3.2** Assume for simplicity that $\mathrm{Im}\varrho$ is the canonical embedding of $\mathrm{PGL}_2(\mathbf{F}_3)$ inside $\mathrm{PGL}_2(\overline{\mathbf{F}}_3)$. The following facts are then easy to deduce from what has been said before in this section:

(a) Any lifting of $\varrho$ into $\mathrm{GL}_2(\mathbf{F}_9)$ has image exactly one of the five groups in the statement of the previous proposition.

(b) The existence of a lifting $\rho$ of $\varrho$ with image one of the five groups above amounts to the solvability of the embedding problem

$$G \longrightarrow \mathcal{S}_4 \simeq \mathrm{Gal}(K/\mathbf{Q})$$

for $G = 2\mathcal{S}_4^+$, $4\mathcal{S}_4^+$, $2\mathcal{S}_4^-$, $4\mathcal{S}_4^-$, $8\mathcal{S}_4^-$, respectively. In this case, any other liftings of $\varrho$ with the same index and the same order of determinant as $\rho$ (see remark 3.1) have also the same image as $\rho$.

(c) Let $\eta \colon G_{\mathbf{Q}} \to \mathrm{GL}_2(\mathbf{F}_9)$ be an octahedral representation with $K_{\overline{\eta}} = K$. Since $K_\eta$ is a solution to someone of the embedding problems in (b), there exists a lifting $\rho$ of $\varrho$ into $\mathrm{GL}_2(\mathbf{F}_9)$ with the same fixed field as $\eta$. Then $\eta$ is conjugate to $\rho$ or to a twist $\rho \otimes \varepsilon$; more precisely, the outer automorphism groups above determine, depending on $\mathrm{Im}\rho$, the possibilities for the character $\varepsilon$:

| $\mathrm{Im}\rho$ | $\varepsilon$ |
|---|---|
| $2\,\mathcal{S}_4^+$ | $\det \rho$ |
| $4\,\mathcal{S}_4^+$ | $\det \rho$, $(\det \rho)^2$, $(\det \rho)^3$ |
| $2\,\mathcal{S}_4^-$ | $\varepsilon_M$ |
| $4\,\mathcal{S}_4^-$ | $\det \rho$, $\varepsilon_M$, $\varepsilon_M \det \rho$ |
| $8\,\mathcal{S}_4^-$ | $(\det \rho)^j$, $\varepsilon_M$, $\varepsilon_M (\det \rho)^j$, $j = 1, 2, 3$ |

All the above twists preserve Im$\rho$, so that $\eta$ is conjugate to a lifting of $\varrho$ into $\mathrm{GL}_2(\mathbf{F}_9)$.

Let us now write down in terms of Hilbert symbols the obstructions to the solvability of the embedding problems defined by the groups in proposition 3.3. If $\varepsilon\colon \mathrm{G}_{\mathbf{Q}} \to F^*$ is the character attached to a quadratic field $\mathbf{Q}(\sqrt{b})$, the corresponding element $[\varepsilon]$ of $\mathrm{Br}_2(\mathbf{Q})$ is given by the Hilbert symbol $(-1, b)$. Provided that it is trivial, all cyclic extensions of $\mathbf{Q}$ of degree 4 containing $\mathbf{Q}(\sqrt{b})$ are then $\mathbf{Q}\left(\sqrt{r(b+x\sqrt{b})}\right)$ for $r, x \in \mathbf{Q}^*$ such that $b - x^2 \in \mathbf{Q}^{*2}$. For a Galois character $\psi$ having fixed field one of these quartics, $[\psi]$ is given by the product $(2, b) \otimes (-1, r)$. From propositions 3.1 and 3.2 we get the following corollary.

**Corollary 3.1** *Consider the embedding problems*

$$G \longrightarrow \mathcal{S}_4 \simeq \mathrm{Gal}(K/\mathbf{Q})$$

*for $G$ one of the groups $2\,\mathcal{S}_4^+$, $4\,\mathcal{S}_4^+$, $2\,\mathcal{S}_4^-$, $4\,\mathcal{S}_4^-$, $8\,\mathcal{S}_4^-$. At least one of them is solvable if and only if there exist $b, r \in \mathbf{Q}^*/\mathbf{Q}^{*2}$ such that*

$$\begin{cases} (-1, b) = 1 \\ w \otimes (-2, d_{K_1}) = (2, b) \otimes (-1, r), \end{cases}$$

*where we use the same notations for $d_{K_1}$ and $w$ as before. If $K$ is not real, the case $G = 2\,\mathcal{S}_4^-$ does never occur, and the following four possibilities for $b$ and $r$ match with the solvability of each one of the remaining embedding problems:*

- $G = 2\,\mathcal{S}_4^+$ *corresponds to* $b = 1$ *and* $r = d_{K_1}$.
- $G = 4\,\mathcal{S}_4^+$ *corresponds to* $b = d_{K_1}$.
- $G = 4\,\mathcal{S}_4^-$ *corresponds to* $b = 1$ *and* $r \neq d_{K_1}$.
- $G = 8\,\mathcal{S}_4^-$ *corresponds to* $b \neq 1, d_{K_1}$.

**Remark 3.3**

(a) The obstructions for $G = 2\,\mathcal{S}_4^+, 4\,\mathcal{S}_4^+$ admit the following reformulation:

- $2\,\mathcal{S}_4^+ \to \mathrm{Gal}(K/\mathbf{Q})$ is solvable if and only if $w = (2, d_{K_1})$.
- $4\,\mathcal{S}_4^+ \to \mathrm{Gal}(K/\mathbf{Q})$ is solvable if and only if $(-1, d_{K_1}) = 1$ and there exists $r \in \mathbf{Q}^*/\mathbf{Q}^{*2}$ such that $w = (-1, r)$.

We must always have $d_{K_1} > 0$ in the second case. If $K$ is not real, then we must have $d_{K_1} < 0$ in the first case, and hence the embedding problems cannot be both solvable. If we replace $\mathbf{Q}$ by any number field $k$, these obstructions remain valid in $\mathrm{Br}_2(k)$ (see [Cre98]).

(b) Given $b_0, r_0 \in \mathbf{Q}^*/\mathbf{Q}^{*2}$ satisfying the conditions of the corollary, any other $b$ (resp. $r$) is obtained from $b_0$ (resp. $r_0$) by multiplying by 2 and primes $p \equiv 1 \pmod{8}$ (resp. by 2 and primes $p \equiv 1 \pmod 4$). In particular, the solvability of $8\,\mathcal{S}_4^- \to \mathrm{Gal}(K/\mathbf{Q})$ is ensured by the existence of any such $b, r \in \mathbf{Q}^*/\mathbf{Q}^{*2}$; we will say that the embedding problem is *solvable with type* $[b, r]$, where $b$ and $r$ are uniquely determined up to the transformations just explained. Notice also that the embedding problem $4\,\mathcal{S}_4^- \to \mathrm{Gal}(K/\mathbf{Q})$ is always solvable whenever $2\,\mathcal{S}_4^+ \to \mathrm{Gal}(K/\mathbf{Q})$ is.

(c) Suppose, as in remark 3.2, that $\mathrm{Im}\varrho$ is the canonical embedding of $\mathrm{PGL}_2(\mathbf{F}_3)$ inside $\mathrm{PGL}_2(\overline{\mathbf{F}}_3)$. Assume also that $\varrho$ is odd. For each pair $b, r \in \mathbf{Q}^*$ (necessarily with $b > 0$ and $r < 0$) satisfying the conditions in the corollary, let $\varepsilon\colon G_{\mathbf{Q}} \to \mathbf{F}_9^*$ be a character with fixed field

$$K_\varepsilon = \begin{cases} \mathbf{Q}(\sqrt{r}) & \text{if } b \in \mathbf{Q}^{*2}, \\ \mathbf{Q}\left(\sqrt{r(b + x\sqrt{b})}\right) & \text{if } b \notin \mathbf{Q}^{*2}, \end{cases}$$

for any $x \in \mathbf{Q}^*$ such that $b - x^2 \in \mathbf{Q}^{*2}$. Then, there is a lifting $\rho\colon G_{\mathbf{Q}} \to \mathrm{GL}_2(\mathbf{F}_9)$ of $\varrho$ with $\det \rho = \varepsilon$. Moreover, $\mathrm{Im}\rho$ is determined by the values of $b$ and $r$ in $\mathbf{Q}^*/\mathbf{Q}^{*2}$:

- If $b = 1$ and $r = d_{K_1}$, then $\mathrm{Im}\rho \simeq 2\,\mathcal{S}_4^+$. This is the only case for which $K_{\det \rho} \subset K$.
- If $b = d_{K_1}$, then $\mathrm{Im}\rho \simeq 4\,\mathcal{S}_4^+$.
- If $b = 1$ and $r \neq d_{K_1}$, then $\mathrm{Im}\rho \simeq 4\,\mathcal{S}_4^-$.
- If $b \neq 1, d_{K_1}$, then $\mathrm{Im}\rho \simeq 8\,\mathcal{S}_4^-$.

In each case, as it is shown in the diagrams below, $K_\rho$ is a quadratic extension of $KK_{\det \rho}$ which solves the corresponding embedding problem attached to $K/\mathbf{Q}$, and every solution to it is of the form $K_\rho$ for some lifting $\rho$ of $\varrho$ into $\mathrm{GL}_2(\mathbf{F}_9)$.

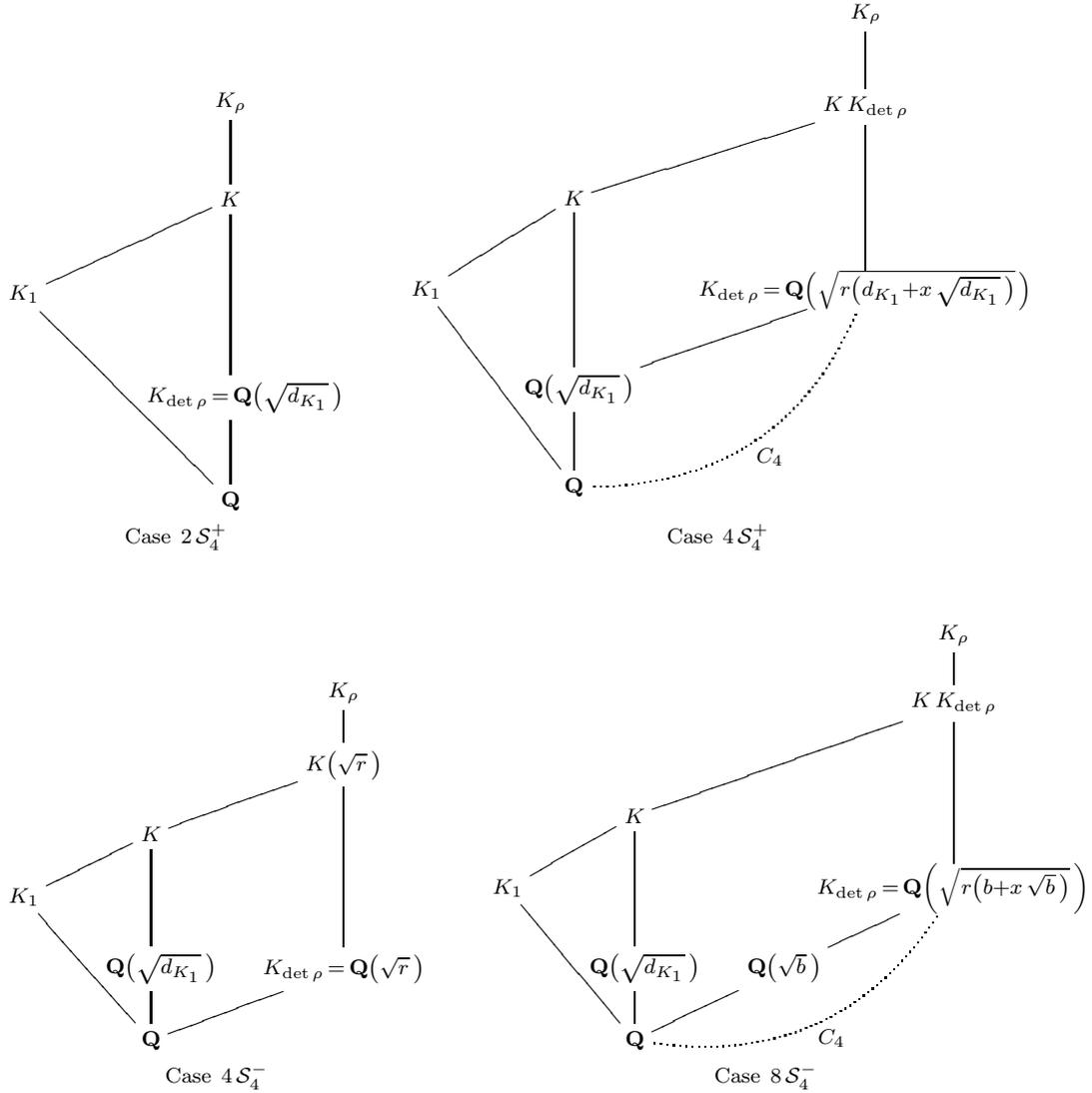

## 4 Principal quartics and embedding problems

We now want to study the solvability of the above embedding problems in the particular case of octahedral extensions defined by principal quartic polynomials.

Before going on, we need to fix some terminology. Let $K_1$ be a quartic extension of a number field $k$ with octahedral normal closure inside $\overline{k}$, and consider the trace linear form $\text{Tr}_{K_1/k}$. Since $K_1/k$ has no quadratic subextensions, the map sending each $\beta \in K_1$ to its minimal polynomial over $k$ becomes a one-to-one correspondence between elements $\beta \in K_1^*$ with $\text{Tr}_{K_1/k}(\beta) = 0$ and *reduced* polynomials $f(X) = X^4 + a\,X^2 + b\,X + c \in k[X]$ defining the extension $K_1/k$. Given such a polynomial $f$, with root $\beta$ in $K_1$, and $m, n, p \in k$ not all three zero, we define the *reduced Tschirnhaus transformation* $\mathrm{Tsc}(f\,;m,n,p\,;X)$ as the minimal polynomial over $k$ of the trace-zero primitive element

$$\gamma = m\,\beta^3 + n\,\beta^2 + p\,\beta + q\,,$$

where the value of $q$ is determined by the condition $\mathrm{Tr}_{K_1/k}(\gamma) = 0$. This reduced polynomial can be computed by means of a resultant:

$$\mathrm{Tsc}(f\,;m,n,p\,;X) = \mathrm{Res}_Y\left(f(Y),\, X - \left(m\,Y^3 + n\,Y^2 + p\,Y + \frac{3\,b\,m + 2\,a\,n}{4}\right)\right).$$

Notice that a reduced quartic polynomial in $k[X]$ defines $K_1/k$ if and only if it can be obtained from $f$ by some Tschirnhaus transformation.

Recall that we say that the quartic $K_1/k$ is *principal* if $f$ admits a Tschirnhaus transformation of the form $g(X) = X^4 + b\,X + c$, for some $b, c \in k$.

**Remark 4.1** If $g(X) = X^4 + b\,X + c$ is a principal polynomial defining $K_1/k$, then for every $m, n \in k$ with $4\,c\,m + 3\,b\,n \neq 0$,

$$\mathrm{Tsc}\left(g\,;m,n,\frac{-3\,b^2\,m^2 + 16\,c\,n^2}{8\,(4\,c\,m + 3\,b\,n)}\,;X\right)$$

is also a principal polynomial defining $K_1/k$.

Let $d_{K_1} \in k^*/k^{*2}$ be the discriminant of the $k$-algebra $K_1$. In $K_1$ we have the trace bilinear form $\mathcal{B}(x,y) = \mathrm{Tr}_{K_1/k}(x\,y)$ and the trace quadratic form

$$\mathcal{Q}(x) = \mathrm{Tr}_{K_1/k}(x^2),$$

which has discriminant $d_{K_1}$. We have already seen how the invariants $\mathrm{disc}\,\mathcal{Q} \in k^*/k^{*2}$ and $w(\mathcal{Q}) \in \mathrm{Br}_2(k)$ attached to this quadratic form appear in the obstruction to the solvability of certain octahedral embedding problems. Let us see how they can be used to characterize principal quartics.

Let $Z = \{x \in K_1 \mid \mathrm{Tr}_{K_1/k}(x) = 0\}$. This is a $k$-vector subspace of $K_1$ whose orthogonal complement with respect to the bilinear form $\mathcal{B}$ is $k$. Since $\mathrm{Tr}_{K_1/k}(1) = 4$, we have that

$$\mathcal{Q} \sim X^2 \oplus \mathcal{Q}_Z\,,$$

where $\mathcal{Q}_Z$ denotes the restriction of $\mathcal{Q}$ to $Z$. Thus, the quadratic forms $\mathcal{Q}$ and $\mathcal{Q}_Z$ have the same discriminant and the same Witt invariant.

For every nonzero $\beta \in Z$ its minimal polynomial $f$ over $k$ is reduced and the coefficient of $X^2$ in $f$ is

$$\frac{\mathrm{Tr}_{K_1/k}(\beta)^2 - \mathrm{Tr}_{K_1/k}(\beta^2)}{2} = -\frac{\mathrm{Tr}_{K_1/k}(\beta^2)}{2} = -\frac{\mathcal{Q}_Z(\beta)}{2}\,.$$

The existence of a polynomial $X^4 + bX + c$ defining $K_1/k$ is equivalent to the existence of an element $\beta \in K_1^*$ such that $\text{Tr}_{K_1/k}(\beta) = \text{Tr}_{K_1/k}(\beta^2) = 0$, namely to the representability of zero by $\mathcal{Q}_Z$, which is given by the element $w(\mathcal{Q}_Z) \otimes (-1, -\text{disc } \mathcal{Q}_Z)$ in $\text{Br}_2(k)$. Putting all together, we have the following result.

**Lemma 4.1** *A quartic extension $K_1/k$ is principal if and only if*

$$w = (-1, -d_{K_1}),$$

*where $w$ is the Witt invariant of the trace quadratic form $\text{Tr}_{K_1/k}(x^2)$ and $d_{K_1}$ is the discriminant of $K_1/k$.*

**Remark 4.2** This condition is explicitly computable from any polynomial $f(X) = X^4 + aX^2 + bX + c$ defining $K_1/k$. The discriminant of $f$ viewed in $k^*/k^{*2}$ equals $d_{K_1}$. If $a \neq 0$ and $d_{K_1} \neq 2a$, the Witt invariant is then $w = \xi \otimes (-1, -d_{K_1})$, where $\xi = (2 a\, d_{K_1},\, 2a^3 + 9b^2 - 8ac)$, so that the condition of the lemma is equivalent to the triviality of the Hilbert symbol $\xi$. Notice also that the rank 3 quadratic form $\mathcal{Q}_Z$ represents 0 if and only if it represents $-d_{K_1}$, and this amounts to the existence of a polynomial $f$ like above such that $d_{K_1} = 2a$.

The following result characterizes principal quartics of $\mathbf{Q}$ with octahedral normal closure in terms of the solvability of the corresponding Galois embedding problem attached to the central extension $8\,\mathcal{S}_4^-$ (see section 3).

**Theorem 4.1** *Let $K_1/\mathbf{Q}$ be a quartic extension with normal closure $K/\mathbf{Q}$ having Galois group isomorphic to $\mathcal{S}_4$. Let $d_{K_1}$ be the squarefree part of the discriminant of $K_1/\mathbf{Q}$ and put*

$$d_{K_1} = \pm 2^\nu d_1\, d_3\, d_5\, d_7,$$

*with $\nu \in \{0,1\}$ and $d_i$ the product of the primes $p \equiv i \pmod 8$ dividing $d_{K_1}$. Then the following conditions are equivalent:*

(a) *The quartic extension $K_1/\mathbf{Q}$ is principal.*

(b) *The embedding problem $8\,\mathcal{S}_4^- \to \text{Gal}(K/\mathbf{Q})$ is solvable with type $[d_5, -d_3]$.*

*Moreover, in this case we have:*

(1) $2\,\mathcal{S}_4^+ \to \text{Gal}(K/\mathbf{Q})$ *is solvable if and only if $d_{K_1} < 0$ and $d_5 = d_7 = 1$.*

(2) $4\,\mathcal{S}_4^+ \to \text{Gal}(K/\mathbf{Q})$ *is solvable if and only if $d_{K_1} > 0$ and $d_3 = d_7 = 1$.*

(3) $4\,\mathcal{S}_4^- \to \text{Gal}(K/\mathbf{Q})$ *is solvable if and only if $d_5 = 1$.*

*Proof.* By lemma 4.1 and corollary 3.1, the obstructions to (a) and (b) are respectively given by the elements $w \otimes (-1, -d_{K_1})$ and $w \otimes (-2, d_{K_1}) \otimes (2, d_5) \otimes (-1, -d_3)$ in $\text{Br}_2(\mathbf{Q})$, where $w$ is the Witt invariant of the trace quadratic form attached to $K_1/\mathbf{Q}$. Recall that, for an odd prime $p$, the Hilbert symbol $(2, p)$ (resp. $(-1, p)$) is trivial if and only if $p \equiv 1, 7 \pmod 8$ (resp. $p \equiv 1, 5 \pmod 8$), while $(2, 2) = (-1, 2) = 1$. Since

$$(-1, -d_{K_1}) \otimes (-2, d_{K_1}) = (-1, -d_3) \otimes (2, d_5),$$

both obstructions are the same. The rest of the theorem is a straightforward consequence of corollary 3.1. □

**Remark 4.3**

(a) If $w = (-1, -d_{K_1})$, then the component at $\infty$ of $[s_4^-] = w \otimes (-2, d_{K_1})$ is $-1$. Hence, for the quartic $K_1/\mathbf{Q}$ to be principal, it cannot be totally real, and then the embedding problem $2\,\mathcal{S}_4^- \to \mathrm{Gal}(K/\mathbf{Q})$ is not solvable.

(b) The condition in (1) is equivalent to the Hilbert symbol $(-2, -3\,d_{K_1})$ being trivial. If $d_{K_1} \neq -3$, this amounts to saying that $-2$ is the norm of some element in the quadratic field $\mathbf{Q}\left(\sqrt{-3\,d_{K_1}}\right)$. Analogously, the condition in (2) is equivalent to $(-1, d_{K_1}) = 1$.

(c) For a discriminant $d_{K_1}$ corresponding to a principal quartic $K_1/\mathbf{Q}$, the graphic below shows, in terms of its factorization in $\mathbf{Q}^*/\mathbf{Q}^{*2}$, which of the embedding problems involved in the theorem are solvable. We have also marked two more sets of discriminants (the reason will become apparent in section 6):

- Those such that $(-1, -3\,d_{K_1}) = 1$, which is equivalent to the existence of some element in $\mathbf{Q}(\sqrt{-3\,d_{K_1}})$ with norm $-1$. In this case, if $d_5 \neq 1$, none of the embedding problems in (1), (2), (3) is solvable.

- Those such that $(2, -3\,d_{K_1}) = 1$, which is equivalent to the existence of some element in $\mathbf{Q}(\sqrt{-3\,d_{K_1}})$ with norm $2$.

Notice that $2$, $-2$ and $-1$ are all norms of elements in $\mathbf{Q}\left(\sqrt{-3\,d_{K_1}}\right)$ if and only if $d_{K_1}$ is of the form $-2^\nu\, 3\, d_1$ in $\mathbf{Q}^*/\mathbf{Q}^{*2}$.

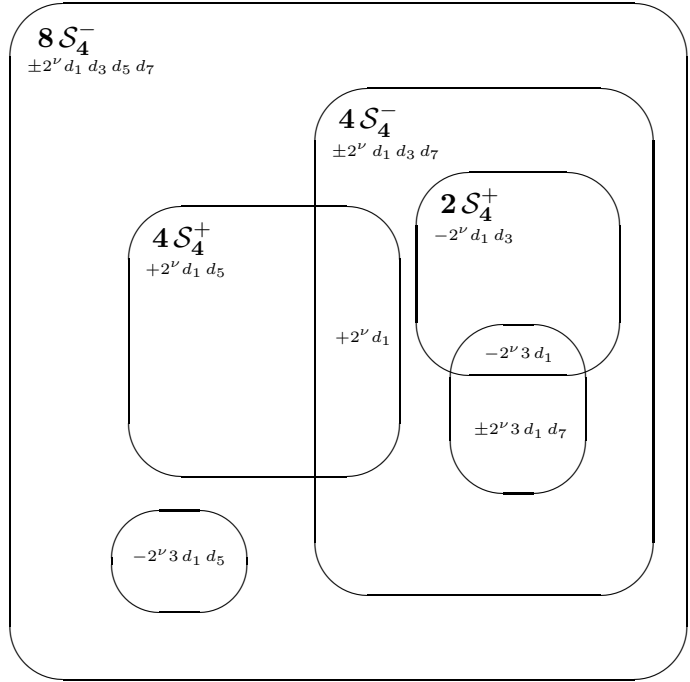

Solvability in $\mathrm{GL}_2(\mathbf{F}_9)$ for odd principal quartics

We finish this section with the following generalization of point (1) in the previous theorem.

**Theorem 4.2** *Let $K$ be an $\mathcal{S}_4$-extension of a number field $k$, and consider the projective representation $\varrho\colon \mathrm{G}_k \to \mathrm{PGL}_2(\mathbf{F}_3)$ obtained from some isomorphism $\mathrm{Gal}(K/k) \simeq \mathrm{PGL}_2(\mathbf{F}_3)$. Let $K_1/k$ be a quartic extension with normal closure $K/k$ and discriminant $d_{K_1} \in k^*/k^{*2}$. Then, any two of the following three conditions imply the third one:*

(a) *The quartic extension $K_1/k$ is principal.*

(b) *The representation $\varrho$ has a lifting to $\mathrm{GL}_2(\mathbf{F}_3)$.*

(c) *The Hilbert symbol $(-2, -3d_{K_1})$ is trivial in $\mathrm{Br}_2(k)$.*

*Proof.* Condition (b) amounts to the solvability of the embedding problem $2\,\mathcal{S}_4^+ \to \mathcal{S}_4 \simeq \mathrm{Gal}(K/k)$. By lemma 4.1 and remark 3.3, the obstructions to (a) and (b) are respectively given by the elements $w \otimes (-1, -d_{K_1})$ and $w \otimes (2, d_{K_1})$ in $\mathrm{Br}_2(k)$, where $w$ is the Witt invariant of the trace quadratic form attached to $K_1/k$. Since $(-1, -d_{K_1}) \otimes (2, d_{K_1}) = (-2, -3\,d_{K_1})$, the result follows. □

**Remark 4.4** For a representation $\varrho$ like in the theorem, $\det \varrho$ is the character $\mathrm{G}_k \to \mathbf{F}_3^*$ attached to the quadratic extension $k\left(\sqrt{d_{K_1}}\right)/k$. It equals the mod 3 cyclotomic character if and only if $d_{K_1} = -3$. In this case, conditions (a) and (b) become equivalent, and the results of [LR95] show that they are also equivalent to the existence of (infinitely many) elliptic curves $E$ defined over $k$ such that $K$ is the fixed field of $\overline{\rho}_{E,3}$, where $\rho_{E,3}$ is the Galois representation attached to the 3-torsion points of $E$. Every lifting of $\varrho$ to $\mathrm{GL}_2(\mathbf{F}_3)$ is then of the form $\rho_{E,3}$ for some elliptic curve $E$. We will give a generalization of that result for the non-cyclotomic case in theorem 6.1 below.

## 5 Q-curves of degree 2 and principal quartics

For every squarefree positive integer $N$, let $X^*(N)$ be the quotient of the modular curve $X_0(N)$ by the automorphism group of Atkin-Lehner involutions. Let $r$ be the number of prime factors of $N$. The curve $X^*(N)$ is defined over $\mathbf{Q}$, and every non-cusp non-CM $\mathbf{Q}$-rational point lifts to a Galois stable set of $2^r$ points on $X_0(N)$ whose $j$-invariants correspond thus to $\mathbf{Q}$-curves having an isogeny of degree dividing $N$ to each one of its conjugate curves. After Elkies's preprint [Elk93] (see also [BL99]), we know that every $\mathbf{Q}$-curve is isogenous to one attached to a rational point of $X^*(N)$ for some $N$. In [GL98], an explicit method is given to parametrize all the $\mathbf{Q}$-curves arising from the rational points of $X^*(N)$ whenever this curve has genus zero or one. Conjecturally, these exhaust all the isogeny classes of $\mathbf{Q}$-curves up to a finite (non-empty) set.

We now want to consider $\mathbf{Q}$-curves of degree 2. The isomorphism classes of such curves are parametrized by the non-cusp non-CM $\mathbf{Q}$-rational points of the genus-zero curve $X^*(2)$ whose preimages on $X_0(2)$ have non-rational $j$-invariant. Recall that the modular function

$$G(z) = \left(\frac{\eta(z)}{\eta(2z)}\right)^{24},$$

$\eta(z)$ being the Dedekind function, is a Hauptmodul on $X_0(2)$ defined over $\mathbf{Q}$. By applying the ideas of [GL98] to the Hauptmodul

$$t = \frac{(G-64)^2}{(G+64)^2}$$

on $X^*(2)$, one obtains the following parametrization of the $j$-invariants associated with the non-cusp $\mathbf{Q}$-rational points on $X^*(2)$:

$$j_t = \frac{64\,(27\,t^2 + 360\,t + 125)}{(t-1)^2} + \frac{128\,(81\,t + 175)}{(t-1)^2}\sqrt{t} = \frac{64\,(5 + 3\sqrt{t}\,)^3}{(-1+\sqrt{t}\,)^2\,(1+\sqrt{t}\,)},$$

for $t$ in $\mathbf{Q} \setminus \{1\}$ and with $j_\infty = 1728$. A Weierstrass model for $j_t$ is

$$C_t \;:\; Y^2 = X^3 - 6\,(5 + 3\,\sqrt{t}\,)\,X + 8\,(7 + 9\,\sqrt{t}\,).$$

Note that there are infinitely many (isomorphism classes of) $\mathbf{Q}$-curves of degree 2 defined over any fixed quadratic field.

We can now show the link between **Q**-curves of degree 2 and principal quartic extensions of **Q**. A non-CM elliptic curve $C$ is a **Q**-curve of degree 2 if and only if it is isomorphic to $C_t$ (changing, if necessary, $C_t$ by its conjugate curve) for some non-square rational value of $t$. In this case, the quadratic field of definition of $C$ is $\mathbf{Q}(\sqrt{t})$, and any Weierstrass model $Y^2 = X^3 + AX + B$ for $C$ over $\mathbf{Q}(\sqrt{t})$ must be a quadratic twist of $C_t$. In particular, the fixed field of $\bar{\rho}_{C,3}$ is the splitting field over $\mathbf{Q}(\sqrt{t})$ of the polynomial

$$f_t(X) = X^4 - 12\,(5 + 3\sqrt{t})\,X^2 + 32\,(7 + 9\sqrt{t})\,X - 12\,(25 + 9t + 30\sqrt{t})$$

whose roots are the $x$-coordinates of the 3-torsion points of $C_t$. Assume that this polynomial has octahedral Galois closure, which amounts to saying that $\rho_{C,3}$ is surjective. Let then $K_C$ be the only $\mathcal{S}_4$-extension of **Q** whose compositum with $\mathbf{Q}(\sqrt{t})$ is the fixed field of $\bar{\rho}_{C,3}$ (see section 2), and put $d_C$ to denote its discriminant. The Tschirnhaus transformation

$$\operatorname{Tsc}\left(f_t\,;\,\frac{1}{72},\,-\frac{1}{36},\,\frac{-37-\sqrt{t}}{36}\,;\,X\right)$$

takes $f_t$ to the principal quartic polynomial

$$g_t(X) \;=\; X^4 + 4\,(t-1)^2\,X - 3\,(t-1)^3\,.$$

Since $g_t$ has rational coefficients, $K_C$ is the splitting field of $g_t$ over **Q**. Moreover, the discriminant of $g_t$ is $-2^8\,3^3\,t\,(t-1)^8$. Hence, $d_C$ is not $-3$ in $\mathbf{Q}^*/\mathbf{Q}^{*2}$, which means that $K_C$ does not contain the cubic roots of unity. Notice also that the quadratic field of definition of $C$ can then be written as $\mathbf{Q}(\sqrt{-3\,d_C})$. We have thus proved half of the following result.

**Proposition 5.1** *For a **Q**-curve $C$ of degree 2 such that $\rho_{C,3}$ is surjective, the octahedral extension $K_C/\mathbf{Q}$ is the Galois closure of a principal quartic with discriminant different from $-3$ in $\mathbf{Q}^*/\mathbf{Q}^{*2}$. Conversely, if we let $K_1/\mathbf{Q}$ be a principal quartic extension with normal closure $K/\mathbf{Q}$ having Galois group isomorphic to $\mathcal{S}_4$ and not containing the cubic roots of unity, there exist infinitely many **Q**-curves $C$ of degree 2 defined over $k = \mathbf{Q}\left(\sqrt{-3\,d_{K_1}}\right)$ such that $K$ equals $K_C$.*

*Proof.* We must only show the second part of the statement. Suppose that $K_1/\mathbf{Q}$ admits a defining polynomial $g(X) = X^4 + bX + c$. Since $\operatorname{disc}(g) = d_{K_1} \neq -3$ in $\mathbf{Q}^*/\mathbf{Q}^{*2}$, the rational value

$$t \;=\; \frac{-\operatorname{disc}(g)}{27\,b^4} \;=\; 1 - \frac{256\,c^3}{27\,b^4}$$

is not a square. Moreover, $k = \mathbf{Q}(\sqrt{t})$ and

$$\operatorname{Tsc}\left(g\,;\,\frac{-4\,(t+3\sqrt{t})}{b\,t},\,\frac{16\,c\,\sqrt{t}}{b^2\,t},\,\frac{-64\,c^2\,\sqrt{t}}{3\,b^3\,t}\,;\,X\right) \;=\; f_t(X)\,.$$

This shows that the fixed field of $\bar{\rho}_{C_t,3}$ is $Kk$; in particular, $C_t$ cannot have complex multiplication. Notice also that $t$ is determined by the ratio $c^3/b^4$, so that any other principal polynomial defining the extension $K_1/\mathbf{Q}$ and giving rise by the above procedure to the same $\mathbf{Q}$-curve $C_t$ must be of the form

$$X^4 + b\,r^3\,X + c\,r^4 = \mathrm{Tsc}(g;0,0,r;X)$$

for some $r \in \mathbf{Q}^*$. Thus, we have infinitely many (isomorphism classes of) $\mathbf{Q}$-curves $C$ of degree 2 defined over $k$ such that $K_C = K$. Indeed, for every rational number $s \ne -4c/3b$, the principal quartic polynomial (cf. remark 4.1)

$$X^4 + b_s\,X + c_s = \mathrm{Tsc}\left(g;1,s,\frac{16\,c\,s^2 - 3\,b^2}{8\,(3\,b\,s + 4\,c)};X\right)$$

gives rise as above to a $\mathbf{Q}$-curve of degree 2, with $j$-invariant

$$\frac{27\,(27\,b_s^8 - 207\,b_s^4\,c_s^3 + 128\,c_s^6)}{2\,c_s^6} + \frac{81\,b_s^2\,(b_s^4 - 3\,c_s^3)}{2\,c_s^6}\sqrt{3\,(27\,b_s^4 - 256\,c_s^3)}\,,$$

satisfying the required conditions. $\square$

**Remark 5.1** The above relation between $\mathbf{Q}$-curves of degree 2 and principal quartics comes naturally from the *Weil restriction process* applied to the 3-torsion polynomial of an elliptic curve. More precisely, fix a $\mathbf{Q}$-curve $C_t$ as before and write

$$f_t\left(X + Y\sqrt{t}\right) = P_t(X,Y) + R_t(X,Y)\sqrt{t}\,,$$

with $P_t, R_t \in \mathbf{Q}[X,Y]$. The resultant

$$\mathrm{Res}_Y\left(P_t(X,Y), R_t(X,Y)\right) \in \mathbf{Q}[X]$$

is a polynomial of degree 16 having an irreducible factor

$$h_t(X) = X^4 - 6\,X^2 + 8\,X + 3\,(8 - 9\,t)$$

with discriminant $\mathrm{disc}(h_t) = -2^8\,3^9\,t\,(t-1)^2$ and Witt invariant

$$w(h_t) = (-1,3) \otimes (t, t-1) = (-1, 3t) = (-1, -\mathrm{disc}(h_t))\,,$$

so that $h_t$ defines a principal quartic extension of $\mathbf{Q}$ (see lemma 4.1). In fact, in the surjective case, the splitting field of $h_t$ is the $\mathcal{S}_4$-extension $K_{C_t}$. The Tschirnhaus transformation

$$\mathrm{Tsc}\left(h_t; -\frac{1}{9}, -\frac{1}{9}, \frac{5}{9}; X\right)$$

takes $h_t$ to the above polynomial $g_t(X) = X^4 + 4\,(t-1)^2\,X - 3\,(t-1)^3$.

# 6 Abelian varieties of $\mathrm{GL}_2$-type attached to principal quartics

In this section we culminate our study of principal octahedral Galois representations characterizing them as the ones arising from **Q**-curves of degree 2. We also describe the **Q**-endomorphism algebras of the abelian varieties of $\mathrm{GL}_2$-type attached to them. Let us first summarize the results of the previous sections.

**Theorem 6.1** *Let $K/\mathbf{Q}$ be an $\mathcal{S}_4$-extension. Fix an isomorphism $\mathrm{Gal}(K/\mathbf{Q}) \simeq \mathrm{PGL}_2(\mathbf{F}_3)$ and consider the projective representation $\varrho$ of $\mathrm{G}_\mathbf{Q}$ obtained from the canonical embedding of $\mathrm{PGL}_2(\mathbf{F}_3)$ inside $\mathrm{PGL}_2(\overline{\mathbf{F}}_3)$,*

$$\varrho \colon \mathrm{G}_\mathbf{Q} \twoheadrightarrow \mathrm{Gal}(K/\mathbf{Q}) \simeq \mathrm{PGL}_2(\mathbf{F}_3) \hookrightarrow \mathrm{PGL}_2(\overline{\mathbf{F}}_3)\,.$$

*Let $K_1/\mathbf{Q}$ be a quartic extension with normal closure $K/\mathbf{Q}$ and discriminant $d_{K_1} \in \mathbf{Q}^*/\mathbf{Q}^{*2}$. Assume that $d_{K_1} \neq -3$ and put $k = \mathbf{Q}\left(\sqrt{-3\,d_{K_1}}\right)$. The following conditions are equivalent:*

(1) *The quartic extension $K_1/\mathbf{Q}$ is principal.*

(2) *The embedding problem $8\mathcal{S}_4^- \to \mathcal{S}_4 \simeq \mathrm{Gal}(K/\mathbf{Q})$ is solvable with type $[d_5, -d_3]$, where $d_i$ is the product of the primes $p \equiv i \pmod{8}$ dividing $d_{K_1}$.*

(3) *For every $b, r \in \mathbf{Q}^*$ with respective squarefree parts $2^{\nu_b} d_5 B$ and $-2^{\nu_r} d_3 R$, where $\nu_b, \nu_r \in \{0,1\}$ and $B$ (resp. $R$) is any product of primes $p \equiv 1 \pmod{8}$ (resp. $p \equiv 1 \pmod{4}$), there exist liftings $\rho \colon \mathrm{G}_\mathbf{Q} \to \mathrm{GL}_2(\mathbf{F}_9)$ of $\varrho$ with*

$$K_{\det \rho} = \begin{cases} \mathbf{Q}(\sqrt{r}) & \text{if } b \in \mathbf{Q}^{*2}, \\ \mathbf{Q}\left(\sqrt{r(b+x\sqrt{b})}\right) & \text{if } b \notin \mathbf{Q}^{*2}, \end{cases}$$

*for any $x \in \mathbf{Q}^*$ such that $b - x^2 \in \mathbf{Q}^{*2}$.*

(4) *There exist (infinitely many) **Q**-curves $C$ of degree 2 defined over $k$ such that $\rho_{C,3}$ is a lifting of the restriction of $\varrho$ to $\mathrm{G}_k$.*

*In this case, the solvability of the embedding problems*

$$G \longrightarrow \mathcal{S}_4 \simeq \mathrm{Gal}(K/\mathbf{Q})\,,$$

*for $G$ one of the groups $2\mathcal{S}_4^+$, $4\mathcal{S}_4^+$, $4\mathcal{S}_4^-$, $8\mathcal{S}_4^-$, depends only on $d_{K_1}$, and all the solutions to them are exactly the fixed fields of the liftings in (3). Moreover, these are all the liftings of $\varrho$ into $\mathrm{GL}_2(\mathbf{F}_9)$.*

To any **Q**-curve $C$ like in (4), we associated in section 2 odd octahedral representations

$$\rho \colon \mathrm{G}_\mathbf{Q} \longrightarrow \mathrm{GL}_2(\mathbf{F}_9)\,,$$

obtained from abelian varieties attached to $C$, and such that $K_{\overline{\rho}} = K$. From remark 3.2, we know that $\rho$ belongs to the set of liftings described in (3). In fact, we will see in 6.2 that they exhaust all such liftings.

Now, we turn to the study of the **Q**-endomorphism algebras for the associated abelian varieties of $\mathrm{GL}_2$-type with minimal dimension. Recall the 2-cocycle class $[c] \in H^2(G_\mathbf{Q}, \mathbf{Q}^*)$ attached in section 2 to any **Q**-curve $C$. In [Que01] and [Que00], Quer computes explicitly its sign component $[c^\pm] \in H^2(G_\mathbf{Q}, \{\pm 1\}) \simeq \mathrm{Br}_2(\mathbf{Q})$ and makes use of it to characterize the fields of complete definition for $C$ up to isogeny. Borrowing his terminology, we will say that a Galois character $\varepsilon \colon G_\mathbf{Q} \to \overline{\mathbf{Q}}^*$ is a *splitting character* for $[c]$ if

$$[\varepsilon] = [c^\pm]$$

in $\mathrm{Br}_2(\mathbf{Q})$. For an abelian number field $L$ containing the minimal field of definition of $C$, the existence of a splitting character for $[c]$ factoring through $\mathrm{Gal}(L/\mathbf{Q})$ amounts to the existence of a quadratic $L$-twist $C'$ of $C$ completely defined over $L$ and such that the Weil restriction $B = \mathrm{Res}_{L/\mathbf{Q}}(C')$ factors over **Q** up to isogeny as a product of pairwise non-isogenous abelian varieties of $\mathrm{GL}_2$-type.

From now on, we assume the notation and the equivalent conditions in theorem 6.1, and we let $C$ be any **Q**-curve of degree 2 defined over $k$. In this case, a Hilbert symbol calculation shows that

$$[c^\pm] = [s_4^-] \otimes [\chi],$$

where $\chi \colon G_\mathbf{Q} \to \{\pm 1\}$ is the mod 3 cyclotomic character and $[s_4^-]$ is the element in $\mathrm{Br}_2(\mathbf{Q})$ attached in section 3 to the $\mathcal{S}_4$-extension $K/\mathbf{Q}$. Hence, from theorem 3.1 we obtain the next result.

**Proposition 6.1** *Let $F/\mathbf{Q}$ be a cyclic extension of degree prime to* 3. *The following conditions are equivalent:*

(i) *Any character $\varepsilon \colon G_\mathbf{Q} \to \overline{\mathbf{Q}}^*$ with fixed field $F$ is a splitting character for the* 2-*cocycle class $[c]$ attached to $C$.*

(ii) *There is a lifting $\rho \colon G_\mathbf{Q} \to \mathrm{GL}_2(\overline{\mathbf{F}}_3)$ of $\varrho$ such that $\chi \det \rho$ has fixed field $F$.*

In particular, the liftings $\rho$ in theorem 6.1 provide us with splitting characters $\varepsilon \colon G_\mathbf{Q} \to \overline{\mathbf{Q}}^*$ for $[c]$: for every $b, r, x$ like in (3), we have such a character with fixed field

$$K_\varepsilon = \begin{cases} \mathbf{Q}(\sqrt{-3\,r}) & \text{if } b \in \mathbf{Q}^{*2}, \\ \mathbf{Q}\left(\sqrt{-3\,r(b + x\sqrt{b}\,)}\right) & \text{if } b \notin \mathbf{Q}^{*2}. \end{cases}$$

The **Q**-curve $C$ has then a twist $C'$ completely defined over $L = k\,K_\varepsilon$ and $B = \mathrm{Res}_{L/\mathbf{Q}}(C')$ is a product of abelian varieties attached to $C$. Following the graphic in remark 4.3, we can use the results of [Que00] to give in each case a splitting character $\varepsilon$ providing minimal **Q**-endomorphism algebras:

(a) If $(-2, -3\,d_{K_1}) = 1$, we can take $b = 1$ and $r = d_{K_1}$, so that $K_\varepsilon = k$. Then $B$ is of $\mathrm{GL}_2$-type and $\mathbf{Q} \otimes \mathrm{End}_{\mathbf{Q}}(B) \simeq \mathbf{Q}(\sqrt{-2}\,)$.

(b) If $(2, -3\,d_{K_1}) = 1$, we can take $b = 1$ and $r = -3$, so that $K_\varepsilon = \mathbf{Q}$. Then $B$ is of $\mathrm{GL}_2$-type and $\mathbf{Q} \otimes \mathrm{End}_{\mathbf{Q}}(B) \simeq \mathbf{Q}(\sqrt{2}\,)$.

(c) If $(-1, -3\,d_{K_1}) = 1$, we can take $b = -3\,d_{K_1}$, so that $K_\varepsilon$ is a quadratic extension of $k$. Then $B$ is $\mathbf{Q}$-isogenous to a product $A_1 \times A_2$, where $A_1$, $A_2$ are abelian varieties of $\mathrm{GL}_2$-type with $\mathbf{Q} \otimes \mathrm{End}_{\mathbf{Q}}(A_j) \simeq \mathbf{Q}(i)$.

(d) If the embedding problem $4\mathcal{S}_4^- \to \mathrm{Gal}(K/\mathbf{Q})$ is solvable, we can take $b = 1$ and $r \neq d_{K_1}, -3$, so that $K_\varepsilon$ is a quadratic different from $k$. Then $B$ is of $\mathrm{GL}_2$-type and $\mathbf{Q} \otimes \mathrm{End}_{\mathbf{Q}}(B) \simeq \mathbf{Q}(\sqrt{2}, \sqrt{-2}\,)$.

(e) If the embedding problem $4\mathcal{S}_4^+ \to \mathrm{Gal}(K/\mathbf{Q})$ is solvable, we can take $b = d_{K_1}$, so that $K_\varepsilon$ is a quadratic extension of $\mathbf{Q}(\sqrt{d_{K_1}}\,)$. Then $B$ factors as in (c), with $\mathbf{Q} \otimes \mathrm{End}_{\mathbf{Q}}(A_j) \simeq \mathbf{Q}(\sqrt{2}, \sqrt{-2}\,)$.

(f) Otherwise, we must always have $b \neq 1, d_{K_1}, -3\,d_{K_1}$, so that $K_\varepsilon$ is a quartic field not containing $k$ or $\mathbf{Q}(\sqrt{d_{K_1}}\,)$. Then $B$ factors again as in (c), with $\mathbf{Q} \otimes \mathrm{End}_{\mathbf{Q}}(A_j) \simeq \mathbf{Q}(\sqrt{2}, \sqrt{-2}\,)$.

We conclude with the following result, which characterizes mod 3 octahedral Galois representations arising from $\mathbf{Q}$-curves of degree 2.

**Theorem 6.2** *Let $\rho\colon \mathrm{Gal}(\overline{\mathbf{Q}}/\mathbf{Q}) \to \mathrm{GL}_2(\mathbf{F}_9)$ be an octahedral representation such that $K_{\overline{\rho}}$ does not contain the cubic roots of unity. Then, the following statements are equivalent:*

(i) *$\rho$ arises from a $\mathbf{Q}$-curve of degree 2.*

(ii) *$K_{\overline{\rho}}$ is the Galois closure of a principal quartic $K_1/\mathbf{Q}$.*

(iii) *The embedding problem $8\,\mathcal{S}_4^- \to \mathrm{Gal}(K_{\overline{\rho}}/\mathbf{Q})$ is solvable with type $[d_5, -d_3]$, where $d_5$ and $d_3$ are defined as above.*

*In this case, $\rho$ arises from infinitely many $\mathbf{Q}$-curves of degree 2 defined over the same quadratic $k = \mathbf{Q}(\sqrt{-3\,d_{K_1}}\,)$. The representation $\overline{\rho}$ has attached abelian varieties of dimension at most 4 and there is always one with $\mathbf{Q}$-endomorphism algebra $\mathbf{Q}(\sqrt{2}, \sqrt{-2}\,)$. Moreover, the following holds for $d = -2, -1, 2$: if there exist elements of norm $d$ in $k$, then $\overline{\rho}$ has attached abelian varieties with $\mathbf{Q}$-endomorphism algebra $\mathbf{Q}(\sqrt{d}\,)$.*

*Proof.* As a consequence of theorem 6.1, we only need to prove that (ii) implies (i). According to section 2 and proposition 5.1, if $K_{\overline{\rho}}$ is principal then there exists a $\mathbf{Q}$-curve $C$ of degree 2 such $K_{\overline{\rho}} = K_{\overline{\rho}_{A,\mathfrak{p}}}$, where $A$ is an abelian variety of $\mathrm{GL}_2$-type attached to a splitting map $\alpha$ for $C$ and with $\mathrm{Im}(\rho_{A,\mathfrak{p}})$ a subgroup of $\mathrm{GL}_2(\mathbf{F}_9)$. Since $\mathrm{PGL}_2(\mathbf{F}_9)$ contains a unique subgroup isomorphic to $\mathcal{S}_4$ up to conjugation, we have that $\rho$ is a twist of $\rho_{A,\mathfrak{p}}$; i.e. $\rho \simeq \rho_{A,\mathfrak{p}} \otimes \widetilde{\varepsilon}$, for some character $\varepsilon\colon \mathrm{G}_{\mathbf{Q}} \to \overline{\mathbf{Q}}^*$. The result then follows by noticing that this twist corresponds to the abelian variety of $\mathrm{GL}_2$-type attached to the splitting map $\alpha\,\varepsilon^{-1}$. □

# 7 Numerical examples

We present four sorts of non-cyclotomic examples. The source for quartic polynomials has been taken from the electronic archives at the web site ftp://megrez.math.u-bordeaux.fr/pub/numberfields/degree4long.

In each one of the following four tables we consider twenty principal quartic extensions $K_1/\mathbf{Q}$ with octahedral normal closure $K/\mathbf{Q}$ ordered by the absolute value of the field discriminants. The first column is labeled $d_{K_1}$ denoting the discriminant of the corresponding quartic field. The second column gives the polynomial as encountered in the Bordeaux data base. Finally, the third column displays a principal quartic polynomial which is Tschirnhaus-equivalent to the source polynomial.

The criteria used to build the tables agree with the diagram in section 4 and are as follows:

- Table 1: The embedding problem $2\,\mathcal{S}_4^+ \to \mathcal{S}_4 \simeq \mathrm{Gal}(K/\mathbf{Q})$ is solvable. To these extensions correspond abelian surfaces with $\mathbf{Q}$-endomorphism algebra $\mathbf{Q}(\sqrt{-2})$. The discriminants marked with the symbol $*$ satisfy $(2, -3\,d_{K_1}) = 1$, so that they also give rise to abelian surfaces with $\mathbf{Q}$-endomorphism algebras $\mathbf{Q}(\sqrt{2})$ and $\mathbf{Q}(i)$.

- Table 2: The embedding problem $4\,\mathcal{S}_4^- \to \mathcal{S}_4 \simeq \mathrm{Gal}(K/\mathbf{Q})$ is solvable but $2\,\mathcal{S}_4^+ \to \mathcal{S}_4 \simeq \mathrm{Gal}(K/\mathbf{Q})$ is not. To these extensions correspond abelian varieties with $\mathbf{Q}$-endomorphism algebra $\mathbf{Q}(\sqrt{-2},\sqrt{2})$. The discriminants marked with the symbol $*$ satisfy $(2, -3\,d_{K_1}) = 1$, so that they also provide abelian surfaces with $\mathbf{Q}$-endomorphism algebra $\mathbf{Q}(\sqrt{2})$.

- Table 3: The embedding problem $4\,\mathcal{S}_4^+ \to \mathcal{S}_4 \simeq \mathrm{Gal}(K/\mathbf{Q})$ is solvable but $4\,\mathcal{S}_4^- \to \mathcal{S}_4 \simeq \mathrm{Gal}(K/\mathbf{Q})$ is not. To these quartics we can attach abelian varieties of $\mathrm{GL}_2$-type with $\mathbf{Q}$-endomorphism algebra $\mathbf{Q}(\sqrt{-2},\sqrt{2})$.

- Table 4: None of the above embedding problems is solvable. In this range of discriminants, none of them satisfies $(-1, -3\,d_{K_1}) = 1$. Again, these quartics give rise to abelian varieties of $\mathrm{GL}_2$-type with $\mathbf{Q}$-endomorphism algebra $\mathbf{Q}(\sqrt{-2},\sqrt{2})$ (see Table 5 below).

| $d_{K_1}$ | $f(x)$ | $f_*(x)$ |
|---|---|---|
| $-283$ | $-1 - x + x^4$ | $-1 + x + x^4$ |
| $-331$ | $-1 + x + x^2 - x^3 + x^4$ | $-43 + 37\,x + x^4$ |
| $-491$ | $-1 + 3\,x - x^2 - x^3 + x^4$ | $-47 + 29\,x + x^4$ |
| $-563$ | $-1 - x + x^2 - x^3 + x^4$ | $-41 + 11\,x + x^4$ |
| $-643$ | $1 - 2\,x - x^3 + x^4$ | $3395 + 4987\,x + x^4$ |
| $-688$ | $-1 - 2\,x + x^4$ | $-1 + 2\,x + x^4$ |
| $-1107\,^*$ | $-1 - 2\,x - x^3 + x^4$ | $3 + 5\,x + x^4$ |
| $-1732$ | $-1 + 3\,x - x^3 + x^4$ | $-1 + 8\,x + x^4$ |
| $-1931$ | $1 - 3\,x + x^4$ | $1 + 3\,x + x^4$ |
| $-2092$ | $-2 - 3\,x + x^2 - x^3 + x^4$ | $-247 + 10\,x + x^4$ |
| $-2096$ | $2 - 2\,x - 2\,x^2 + x^4$ | $343 + 842\,x + x^4$ |
| $-2116$ | $-2 + x^2 - x^3 + x^4$ | $-529 + 184\,x + x^4$ |
| $-2243$ | $-1 - 3\,x - x^2 - x^3 + x^4$ | $-15 + 13\,x + x^4$ |
| $-2284$ | $-4 + 2\,x + 2\,x^2 - 2\,x^3 + x^4$ | $4736 + 1700\,x + x^4$ |
| $-2563$ | $-1 - x + 4\,x^2 - 2\,x^3 + x^4$ | $59 + 109\,x + x^4$ |
| $-2608$ | $-2 - 2\,x - 2\,x^2 + x^4$ | $-220 + 148\,x + x^4$ |
| $-2619\,^*$ | $3 + 3\,x - 3\,x^2 - x^3 + x^4$ | $27 + 23\,x + x^4$ |
| $-2764$ | $-2 + x - 3\,x^2 - x^3 + x^4$ | $7256 + 1390\,x + x^4$ |
| $-2787\,^*$ | $-3 - 5\,x - x^2 - x^3 + x^4$ | $-3963 + 1325\,x + x^4$ |
| $-2843$ | $-1 + 2\,x + 2\,x^2 - x^3 + x^4$ | $-10151 + 739\,x + x^4$ |

Table 1
$2\,\mathcal{S}_4^+ \to \mathcal{S}_4 \simeq \mathrm{Gal}(K/\mathbf{Q})$

| $d_{K_1}$ | $f(x)$ | $f_*(x)$ |
|---:|---:|---:|
| 892 | $2 - x^2 - x^3 + x^4$ | $79 + 8x + x^4$ |
| 1436 | $2 + 3x^2 - x^3 + x^4$ | $407 + 32x + x^4$ |
| $-1472$ | $-2 + 2x^2 - 2x^3 + x^4$ | $4 + 8x + x^4$ |
| $-1984$ | $-1 + 2x + 2x^2 - 2x^3 + x^4$ | $-4 + 8x + x^4$ |
| 2021 | $2 - x + x^4$ | $2 + x + x^4$ |
| 2429 | $3 + 2x - 2x^2 - x^3 + x^4$ | $327 + 31x + x^4$ |
| $-2488$ | $-2 + 2x + x^2 - x^3 + x^4$ | $16139 + 2516x + x^4$ |
| $3024\,^*$ | $3 - 2x^3 + x^4$ | $12 + 4x + x^4$ |
| $3261\,^*$ | $3 - 2x + 2x^2 - x^3 + x^4$ | $6297 + 1163x + x^4$ |
| $-4648$ | $-4 - 6x + x^2 + x^4$ | $419 + 388x + x^4$ |
| 4892 | $3 - x - x^3 + x^4$ | $54039 + 2336x + x^4$ |
| 5056 | $3 - 2x + 5x^2 - 2x^3 + x^4$ | $7 + 4x + x^4$ |
| 5056 | $2 + 4x + x^2 - 2x^3 + x^4$ | $156 + 32x + x^4$ |
| $5076\,^*$ | $4 - x + 3x^2 - x^3 + x^4$ | $1512 + 108x + x^4$ |
| 5216 | $2 - 4x + 3x^2 + x^4$ | $11 + 4x + x^4$ |
| $-5348$ | $-1 - x + 4x^2 - x^3 + x^4$ | $-47681 + 424x + x^4$ |
| $5373\,^*$ | $5 + x - 3x^2 - x^3 + x^4$ | $483 + 37x + x^4$ |
| $-5432$ | $-3 + x - x^3 + x^4$ | $-16517 + 2020x + x^4$ |
| $-5816$ | $-1 + 5x - 2x^2 - x^3 + x^4$ | $1227 + 364x + x^4$ |
| $5853\,^*$ | $3 + x - x^2 - x^3 + x^4$ | $11949 + 1843x + x^4$ |

Table 2
$4\mathcal{S}_4^- \to \mathcal{S}_4 \simeq \mathrm{Gal}(K/\mathbf{Q})$

| $d_{K_1}$ | $f(x)$ | $f_*(x)$ |
|---|---|---|
| 229 | $1 - x + x^4$ | $1 + x + x^4$ |
| 592 | $1 - 2x + 2x^2 + x^4$ | $4 + 4x + x^4$ |
| 788 | $2 - 2x + 2x^2 - x^3 + x^4$ | $1184 + 326x + x^4$ |
| 1076 | $2 - 3x + 3x^2 - x^3 + x^4$ | $2237 + 566x + x^4$ |
| 1229 | $3 - x + 3x^2 - x^3 + x^4$ | $2567 + 603x + x^4$ |
| 1396 | $2 - x + x^2 - x^3 + x^4$ | $73 + 46x + x^4$ |
| 1492 | $2 + 2x^2 - x^3 + x^4$ | $13 + 10x + x^4$ |
| 1556 | $2 + x - x^2 - x^3 + x^4$ | $1352 + 172x + x^4$ |
| 1616 | $2 - 2x + x^4$ | $2 + 2x + x^4$ |
| 1940 | $2 - x^3 + x^4$ | $8 + 4x + x^4$ |
| 2213 | $3 - x - 2x^2 + x^4$ | $4589 + 669x + x^4$ |
| 2349 | $3 - 3x + 3x^2 - x^3 + x^4$ | $21 + 11x + x^4$ |
| 2597 | $1 - x + 4x^2 + x^4$ | $1715 + 427x + x^4$ |
| 2677 | $3 - 2x - x^3 + x^4$ | $643 + 185x + x^4$ |
| 2836 | $2 - 4x + 4x^2 - x^3 + x^4$ | $16 + 14x + x^4$ |
| 2917 | $2 - x + 4x^2 + x^4$ | $2784 + 247x + x^4$ |
| 3028 | $2 + 2x - x^3 + x^4$ | $112 + 30x + x^4$ |
| 3221 | $2 - 3x + 2x^2 + x^4$ | $4766 + 137x + x^4$ |
| 3229 | $3 - 4x + 2x^2 - x^3 + x^4$ | $21 + 5x + x^4$ |
| 3316 | $4 - 2x + 2x^2 + x^4$ | $2361 + 526x + x^4$ |

Table 3
$4\mathcal{S}_4^+ \to \mathcal{S}_4 \simeq \mathrm{Gal}(K/\mathbf{Q})$

| $d_{K_1}$ | $f(x)$ | $f_*(x)$ |
|---|---|---|
| $-848$ | $1 - 2x - x^2 + x^4$ | $19 + 16x + x^4$ |
| $-976$ | $-1 + 3x^2 - 2x^3 + x^4$ | $-13 + 32x + x^4$ |
| $-1192$ | $-1 + x + 2x^2 - x^3 + x^4$ | $-13 + 20x + x^4$ |
| $-1456$ | $1 - 2x - 2x^2 + x^4$ | $-4 + 4x + x^4$ |
| $-1856$ | $-1 - 2x + x^2 - 2x^3 + x^4$ | $-404 + 64x + x^4$ |
| $-1963$ | $-1 - 2x + 2x^2 - x^3 + x^4$ | $-25 + 617x + x^4$ |
| $-2051$ | $1 + 3x - x^2 - x^3 + x^4$ | $1837 + 571x + x^4$ |
| $-2219$ | $-1 - x + 3x^2 - x^3 + x^4$ | $12571 + 2097x + x^4$ |
| $-2443$ | $-1 - 3x + x^4$ | $-1 + 3x + x^4$ |
| $-2480$ | $-2 - 2x + x^4$ | $-2 + 2x + x^4$ |
| $-3052$ | $-2 - x + x^2 - x^3 + x^4$ | $-211 + 38x + x^4$ |
| $3760$ | $5 - x^2 - 2x^3 + x^4$ | $1083 + 328x + x^4$ |
| $-4432$ | $3 - 4x - x^2 - 2x^3 + x^4$ | $-13 + 120x + x^4$ |
| $-4595$ | $1 - 3x - 2x^3 + x^4$ | $2713 + 5743x + x^4$ |
| $-4780$ | $2 - 3x - 3x^2 - x^3 + x^4$ | $293 + 1046x + x^4$ |
| $-4907$ | $-1 - 4x - 2x^2 - x^3 + x^4$ | $-11 + 17x + x^4$ |
| $-5224$ | $-1 - 5x + 2x^2 - x^3 + x^4$ | $227 + 196x + x^4$ |
| $-5260$ | $-2 + 3x + x^2 - x^3 + x^4$ | $-112 + 44x + x^4$ |
| $5724$ | $6 - 3x^2 - x^3 + x^4$ | $10407 + 1328x + x^4$ |
| $-5755$ | $1 + 4x - x^3 + x^4$ | $-1807 + 37x + x^4$ |

Table 4

$8\,\mathcal{S}_4^- \to \mathcal{S}_4 \simeq \mathrm{Gal}(K/\mathbf{Q})$

In Table 5, we give five examples of principal quartic extensions $K_1/\mathbf{Q}$ with discriminant $d_{K_1}$ satisfying $(-1, -3\,d_{K_1}) = 1$ and such that none of the embedding problems corresponding to Tables 1, 2 and 3 is solvable. According to the notation in section 4, this amounts to saying that $d_{K_1} = -2^\nu\,3\,d_1\,d_5$, with $d_5 \neq 1$. These quartic fields give rise to abelian varieties of $\mathrm{GL}_2$-type with $\mathbf{Q}$-endomorphism algebra $\mathbf{Q}(i)$.

| $d_{K_1}$ | $f(x)$ | $f_*(x)$ |
|---|---|---|
| $-7155$ | $-3 - x^3 + x^4$ | $-27 + 9x + x^4$ |
| $-8640$ | $-6 + 6x^2 - 2x^3 + x^4$ | $-12 + 8x + x^4$ |
| $-11448$ | $-8 + 8x - 3x^2 - x^3 + x^4$ | $-21 + 12x + x^4$ |
| $-13176$ | $-2 + 2x + 3x^2 - x^3 + x^4$ | $27 + 324x + x^4$ |
| $-14715$ | $-3 - 3x - 3x^2 - x^3 + x^4$ | $-87 + 41x + x^4$ |

Table 5

As a final remark, it is worth to point out that we did not pay attention to the problem of minimizing the coefficients in $f_*(x)$. This problem deserves further investigation since it is related to the problem of choosing a $\mathbf{Q}$-curve of degree 2 with minimal bad reduction. A possible approach would take into account the reducing techniques developed by Cremona [Cre99].

Departament de Matemàtica Aplicada II
Universitat Politècnica de Catalunya
Pau Gargallo 5, E-08028 Barcelona